\documentclass[conference]{IEEEtran}

\usepackage{ifpdf}
\usepackage{url}
\ifpdf
\usepackage[pdftex]{graphicx}
\else
\usepackage[dvips]{graphicx}
\fi
\usepackage{epstopdf}
\DeclareGraphicsExtensions{.eps}
\usepackage[font=footnotesize]{subfig}
\usepackage{caption}
\usepackage{comment}
\usepackage[cmex10]{amsmath}
\usepackage{amsthm,amsfonts,amssymb}
\usepackage[noadjust]{cite}
\usepackage{graphics}
\usepackage{color}
\usepackage[font={small}]{caption}
\usepackage{multirow}
\usepackage{accents}
\newcommand{\ubar}[1]{\underaccent{\bar}{#1}}

\newenvironment{ldescription}[1]
{\begin{list}{}%
		{\renewcommand\makelabel[1]{##1\hfill}%
			\settowidth\labelwidth{\makelabel{#1}}%
			\setlength\leftmargin{\labelwidth}
			\addtolength\leftmargin{\labelsep}}}
	{\end{list}}

\usepackage[utf8]{inputenc}

\begin{document}
	
\title{Impact of Local Transmission Congestion on Energy Storage Arbitrage Opportunities \vspace{-10pt}}
	
\author{	
\IEEEauthorblockN{Yishen~Wang\IEEEauthorrefmark{1},
				  Yury~Dvorkin\IEEEauthorrefmark{2},
				  Ricardo~Fern\'andez-Blanco,
				  Bolun~Xu\IEEEauthorrefmark{1},
				  Daniel S. Kirschen\IEEEauthorrefmark{1}}
\IEEEauthorblockA{\IEEEauthorrefmark{1}Department of Electrical Engineering, University of Washington, Seattle, WA, USA}
\IEEEauthorblockA{\IEEEauthorrefmark{2}Tandon School of Engineering, New York University, New York, NY, USA}
\IEEEauthorblockA{Email: \{ywang11, xubolun, kirschen\}@uw.edu, Ricardo.FCarramolino@gmail.com, dvorkin@nyu.edu} 	\vspace{-10pt}
}

\maketitle
	
\begin{abstract}
Reduced installation and operating costs give energy storage systems an opportunity to participate actively and profitably in electricity markets. In addition to providing ancillary services, energy storage systems can also arbitrage temporal price differences. Congestion in the transmission network often accentuates these price differences and will under certain circumstances enhance the profitability of arbitrage. On the other hand, congestion may also limit the ability of a given storage device to take advantage of arbitrage opportunities. 

This paper analyzes how transmission congestion affects the profitability of arbitrage by storage devices in markets with perfect and imperfect competition. Imperfect competition is modeled using a bilevel optimization where the offers and bids submitted by the storage devices can alter the market outcome. Price-taker and price-maker assumptions are also investigated through market price duration curves. This analysis is based on simulating an entire year of market operation on the IEEE Reliability Test system.
 
\end{abstract}

\begin{IEEEkeywords}
Energy Storage, Bilevel Programming, MPEC, Network-Constrained Market Clearing, Transmission Congestion
\end{IEEEkeywords} 
	
\IEEEpeerreviewmaketitle
\vspace{-10pt}
\section*{Nomenclature}
\subsection{Sets and Indices}
\begin{ldescription}{$xxxxxxxxxx$}
	\item [$B$] Set of buses, indexed by $b$.
	\item [$I$] Set of conventional generators, indexed by $i$.
	\item [$L$] Set of transmission lines, indexed by $l$.
	\item [$H$] Set of storage devices, indexed by $h$.
	\item [$W$] Set of wind generators, indexed by $w$. 
	\item [$T$] Set of time intervals, indexed by $t$. 
	\item [$f(l)/t(l)$] Indices of sending/receiving buses of line $l$.
\end{ldescription}

\subsection{Variables}
\begin{ldescription}{$xxxxxxxxxx$}
	\item [$chs_{h,t} / dis_{h,t}$] Quantity of charging bid/discharging offer of storage $h$ at time $t$, MW.
	\item [$\rho_{h,t}^{chs} / \rho_{h,t}^{dis}$] Price of charging bid/discharging offer of storage $h$ at time $t$, \$/MWh.
	\item [$d_{b,t}$] Load at bus $b$ at time $t$, MW.
	\item [$p_{i,t}^G$] Power output of conventional generator $i$ at time $t$, MW.
	\item [$p_{w,t}^W$] Power output of wind generator $w$ at time $t$, MW. 
	\item [$pf_{l,t}$] Power flow on line $l$ at time $t$, MW. 
	\item [$q_{h,t}^{dis/chs}$] Cleared discharging/charging rate for storage $h$ at time $t$, MW. 
	\item [$SoC_{h,t}$] State of charge of storage $h$ at time $t$, MWh.
	\item [$x_{h,t}^{dis/chs}$] Binary variable representing the discharging/charging status of energy storage $h$ at time $t$. This variable is equal to  1 if discharging/charging and 0 otherwise.
	\item [$\theta_{b,t}$] Voltage phase angle at bus $b$ at time $t$, radians.
	\item [$\lambda_{b,t}$] Dual variables associated with the power balance constraint at bus $b$ at time $t$.
\end{ldescription}	
\vspace{-10pt}
\subsection{Parameters}
\begin{ldescription}{$xxxxxxxxxx$}
	\item [$C_{h}^{bid} $] Charging price bid to reflect storage owners' willingness to pay for charging and marginal charging cost of storage $h$, \$/MW.
	\item [$C_{h}^{dis/chs} $]  Marginal discharging/charging cost of storage $h$, \$/MW.
	\item [$C_{i}^{G} / C_{w}^{W}$] Offer price of conventional generator $i$ / wind generator $w$, \$/MWh.
	\item [$C_{b}^{D} $] Bid price of consumers at bus $b$, \$/MWh.
	\item [$\overline{dis}_{h} / \overline{chs}_{h} $] Discharging/charging rate limits of storage $h$, MW.
	\item [$\eta_{h}^{dis} / \eta_{h}^{chs}  $] Discharging/charging efficiency of storage $h$.
	\item [$\ubar{D}_{b,t} / \bar{D}_{b,t}$] Min/max bounds on consumers bid at bus $b$ at time $t$, MW.
	\item [$\bar{F}_{l}$] Power flow limit on line $l$, MW.
	\item [$\ubar{P}_{i} / \bar{P}_{i}$] Min/max power output limits of generator $i$, MW.
	\item [$\underline{SoC}_{h} / \overline{SoC}_{h} $] Min/max state-of-charge of storage $h$, MWh.
	\item [$SoC_{h}^{init}$] Initial state of charge of storage $h$, MWh.
	\item [$X_{l}$] Reactance of line $l$.
	\item [$WF_{w,t}$] Forecast output power of wind generator $w$ at time $t$, MW.
\end{ldescription}
	\vspace{-10pt}
\section{Introduction}

Due to its declining installation and operation cost, energy storage systems are likely to play an increasing role not only in power system operation but also in electricity markets. Because of its flexibility, storage can be used for energy arbitrage \cite{hrvoje_2015}, transmission congestion relief \cite{vargas_2015}, reserves \cite{pozo_2014_storage_uc}, frequency regulation \cite{low_2016_frequency}, post-contingency corrective actions \cite{wen_2016} and other services. In order to encourage the deployment of energy storage, market rules are being revised to support products and pricing schemes better suited to the technical characteristics and constraints of storage systems \cite{xu_2016_reg}. 

In a vertically-integrated environment, energy storage is used to minimize the system operating cost. For example, Pandzic \emph{et al.} \cite{hrvoje_2015}, Qiu \emph{et al.} \cite{qiu_2016} and Fernández-Blanco \emph{et al.} \cite{ricardo_wecc} investigated the optimal siting and sizing of storage systems used for spatio-temporal arbitrage. The authors of \cite{dvorkin_2016_bilevel} proposed a bilevel storage planning strategy that guarantees the recovery of the investment costs. Wogrin \emph{et al.} \cite{worgin_2016} proposed to allocate storage for load-shifting and regulation services. Wen \emph{et al.} \cite{wen_2016} demonstrated the effectiveness of storage for post-contingency corrective actions. However, minimizing the operating cost, or more generally maximizing the benefits to the system, typically does not lead to a charging/discharging schedule that maximizes the profits of a merchant storage owner. Operation and planning of storage in a competitive market environment is therefore gaining increasing attention. Castillo \emph{et al.} \cite{castillo_2013}, Shafiee \emph{et al.} \cite{shafiee_2016}, and Mohsenian-Rad \cite{rad_2016} proposed techniques to maximize the operating profits of storage in a decentralized environment. Ding \emph{et al.} \cite{ding_2016} coordinated the operation of storage and wind farms in a rolling real-time market. Xu \emph{et al.} \cite{xu_2016_reg} analyzed the batteries' profits for regulation services considering its degradation costs \cite{xu_2016_degrad, ortega_2014}. Khani \emph{et al.} \cite{khani_2016_congestion} proposed to perform arbitrage while relieving transmission congestion. Relative to \cite{shafiee_2016, rad_2016, ding_2016, xu_2016_degrad, ortega_2014, khani_2016_congestion}, this paper models impacts of local and system-wide transmission congestion and  bidding of other market participants.

A large number of papers discuss how network constraints and transmission congestion affect electricity markets and influence the bidding strategy of market participants, e.g.   \cite{wu_2007, tomsovic_2003, gao_2010}. However, previous work focused on system-wide transmission congestion, while this paper examines how storage can affect or be affected by local transmission congestion. This paper makes three contributions: 

\begin{enumerate}
    \item A network-constrained market-clearing mechanism with storage participation is analyzed under perfect and imperfect competition. Perfect competition is modeled as an economic dispatch problem, while imperfect competition is modeled using a bilevel approach where the storage owner  behaves strategically. 
    \item The effect of local transmission congestion on the annual operating profit of storage is quantified to show how it affects the behavior of storage owners.
    \item The importance of optimizing siting and sizing decisions on the  storage profitability  is discussed. 
\end{enumerate}

The rest of the paper is organized as follows. Section II presents the mathematical formulation of a network-constrained market clearing with storage participation under perfect and imperfect competition. Section III presents a detailed case study based on a year-long simulation of storage operation in the IEEE RTS. The influence of local transmission capacity, price-taker/price-maker assumption and storage siting decisions are discussed. Section IV concludes the paper.
	\vspace{-10pt}
\section{Mathematical Formulation}
This section provides the mathematical formulation of a network-constrained market clearing with storage participation under perfect and imperfect competition.
 	\vspace{-10pt}
\subsection{Market Clearing under Perfect Competition}
Under perfect competition, each market participant bids or offers at its marginal cost because none of them is assumed to be able to exercise market power. The ISO then clears the market in a way that maximizes the social welfare. Storage owners participate in the market by  submitting bids and offers that reflect their willingness to charge and discharge. These bids and offers should take into account how battery cycling affects the life of the battery, i.e. its incremental degradation cost. Market clearing takes the form of a network-constrained economic dispatch problem. Storage owners and other market participants pay or are paid based on their locational marginal prices (LMPs). 

	\begin{equation} \label{eq:ed_obj}
	\begin{aligned}
		\max \quad Obj^{SO} & = \sum_{b,t}^{} C_b^{D} \, d_{b,t} + \sum_{h,t}^{} C_h^{bid} \, q_{h,t}^{chs}   \\
							& - \sum_{i,t}^{} C_i^G \, p_{i,t}^G - \sum_{h,t}^{} C_h^{dis} \, q_{h,t}^{dis}
	\end{aligned} 
	\end{equation}
	\begin{equation} \label{eq:ed_chs}
		0 \leq q_{h,t}^{chs} \leq \overline{chs}_h x_{h,t}^{chs}, \, \forall h, \forall t
	\end{equation}
	\begin{equation} \label{eq:ed_dis}
		0 \leq q_{h,t}^{dis} \leq \overline{dis}_h x_{h,t}^{dis}, \, \forall h, \forall t
	\end{equation}
	\begin{equation} \label{eq:ed_dis_state}
		x_{h,t}^{dis} + x_{h,t}^{chs} \leq 1, \, \forall h, \forall t
	\end{equation}
	\begin{equation} \label{eq:ed_binary}
		x_{h,t}^{dis}, \, x_{h,t}^{chs} \in \{0,1\}, \, \forall h, \forall t
	\end{equation}
	\begin{equation} \label{eq:ed_soc_t} 
		SoC_{h,t} = SoC_{h,t-1} + q_{h,t}^{chs} \, \eta_h^{chs} - q_{h,t}^{dis} / \eta_h^{dis}, \, \forall t>1, \forall h
	\end{equation}
	\begin{equation} \label{eq:ed_soc0} 
		SoC_{h,t} = SoC_{h}^{init} + q_{h,t}^{chs} \, \eta_h^{chs} - q_{h,t}^{dis} / \eta_h^{dis}, \, t=1, \forall h
	\end{equation}
	\begin{equation} \label{eq:ed_soc} 
		\underline{SoC}_h \leq SoC_{h,t} \leq \overline{SoC}_h, \, \forall h, \forall t
	\end{equation}
	\begin{equation} \label{eq:ed_socT0}
		SoC_{h,N_T}  = SoC_{h}^{init}, \, \forall h
	\end{equation}
	\begin{equation} \label{eq:ed_gen}
		0  \leq p_{i,t}^G \leq \bar{P}_i, \, \forall i, \forall t
	\end{equation}
	\begin{equation} \label{eq:ed_wind}
		0 \leq p_{w,t}^W \leq WF_{w,t}, \, \forall w, \forall t
	\end{equation}
	\begin{equation} \label{eq:ed_load}
		\ubar{D}_{b,t} \leq d_{b,t} \leq \bar{D}_{b,t}, \, \forall b, \forall t
	\end{equation}
	\begin{equation} \label{eq:ed_powerbalance}
		\begin{aligned}
			& d_{b,t} + \sum_{b | f(l)=b}^{} pf_{l,t} - \sum_{b | t(l)=b}^{} pf_{l,t} = \\ 
			& \sum_{i}^{} p_{b(i),t}^G + \sum_{w} p_{b(w),t}^W + \sum_{h}^{} (q_{b(h),t}^{dis} - q_{b(h),t}^{chs}), \, \forall b,\forall t
		\end{aligned}
	\end{equation}
	\vspace{-10pt}
	\begin{equation} \label{eq:ed_pfdef}
		pf_{l,t} = \frac{1}{X_{l}} (\theta_{f(l),t} - \theta_{t(l),t}), \, \forall l, \forall t
	\end{equation}
	\begin{equation} \label{eq:ed_pf}
		- \bar{F_l} \leq pf_{l,t} \leq \bar{F_l}, \, \forall l, \forall t 
	\end{equation}	
	\begin{equation} \label{eq:ed_angle}
		-\pi \leq \theta_{b,t} \leq \pi, \forall b \neq \text{ref}, \forall t  	
	\end{equation}
	\begin{equation} \label{eq:ed_vref}
		\theta_{b,t} = 0, \, b = \text{ref}, \forall t.	\vspace{-10pt}
	\end{equation}

The objective function \eqref{eq:ed_obj} maximizes the social welfare which includes generation and storage discharging offers as well as the consumers and storage charging  bids. Constraints \eqref{eq:ed_chs}--\eqref{eq:ed_dis} set the limits on the storage bid/offer quantities. Eqs. \eqref{eq:ed_dis_state}--\eqref{eq:ed_binary} prevent simultaneous charging and discharging by enforcing constraints on binary variables. Constraints \eqref{eq:ed_soc_t}--\eqref{eq:ed_soc} track the state of charge (SoC) and enforce the limits on their operating range. Constraint \eqref{eq:ed_socT0} forces the final SoC to be identical to the initial SoC. Constraints \eqref{eq:ed_gen}--\eqref{eq:ed_load} enforce the lower and upper bounds on the thermal generation, wind farms, and consumers. Constraint \eqref{eq:ed_powerbalance} is the nodal power balance. Constraint \eqref{eq:ed_pfdef} calculates the power flows using a dc load flow model. Constraints \eqref{eq:ed_pf}--\eqref{eq:ed_vref} enforce the limits on the line flows and the voltage angles. 

\subsection{Market Clearing under Imperfect Competition}
In an imperfectly competitive market, storage can achieve larger profits through strategic bidding and offering. This strategic behavior can be modeled as a bilevel program that captures the interactions between the merchant storage and the ISO. In the upper level, storage maximizes its operating profits and determines the price/quantity bids and offers to be submitted to the ISO. The lower level represents a network-constrained market clearing as described in the previous section. The accepted bids, offers and market-clearing locational marginal prices are fed back to the upper level where they are used to calculate the profits of storage.   
	\begin{equation} \label{eq:up_obj}
		\max \, Obj^{ESS} = \sum_{h,t} [\lambda_{b(h),t} (q_{h,t}^{dis}-q_{h,t}^{chs}) - C_h^{dis} q_{h,t}^{dis} - C_h^{chs} q_{h,t}^{chs}]
	\end{equation}

	\begin{equation} \label{eq:up_offer}
		\rho_{h,t}^{dis}, \rho_{h,t}^{chs} \geq 0, \, \forall h, \forall t 
	\end{equation}
	\begin{equation} \label{eq:up_chs}
		0 \leq chs_{h,t} \leq \overline{chs}_h x_{h,t}^{chs}, \, \forall h, \forall t
	\end{equation}
	\begin{equation} \label{eq:up_dis}
		0 \leq dis_{h,t} \leq \overline{dis}_h x_{h,t}^{dis}, \, \forall h, \forall t
	\end{equation}
	\begin{equation} \label{eq:up_same_es}
		\text{Constraints \eqref{eq:ed_dis_state}--\eqref{eq:ed_socT0}}
	\end{equation}
	\begin{equation} \label{eq:lo_obj}
	\begin{aligned}
		\lambda_{b,t}, q_{h,t}^{dis}, q_{h,t}^{chs} \in \arg\max \bigg\{ SW = \sum_{b,t}^{} C_b^{D} \, d_{b,t} \\
		+ \sum_{h,t}^{} \rho_{h,t}^{chs} \, q_{h,t}^{chs} - \sum_{h,t}^{} \rho_{h,t}^{dis} \, q_{h,t}^{dis} - \sum_{i,t}^{} C_i^G \, p_{i,t}^G
	\end{aligned} 
	\end{equation}
		\begin{equation} \label{eq:lo_chs}
			0 \leq q_{h,t}^{chs} \leq chs_{h,t}, \, \forall h, \forall t
		\end{equation}
		\begin{equation} \label{eq:lo_dis}
			0 \leq q_{h,t}^{dis} \leq dis_{h,t}, \, \forall h, \forall t
		\end{equation}
		\begin{equation} \label{eq:lo_same}
			\text{Constraints \eqref{eq:ed_gen}--\eqref{eq:ed_vref}} \bigg\}.
		\end{equation}

The upper-level objective function \eqref{eq:up_obj} maximizes the storage profits based on the LMPs and the cleared charging/discharging quantities. Constraint \eqref{eq:up_offer} enforces the non-negativity of bid and offer prices submitted by the storage. This constraint could be relaxed in systems with high renewable penetration where negative electricity prices can occur. Constraints \eqref{eq:up_chs}--\eqref{eq:up_dis} set the limits on the quantity bids and offers. Other constraints on storage operation are identical to those described in the previous section. LMPs and cleared quantities are obtained from the lower-level problem. This lower-level problem maximizes the social welfare \eqref{eq:lo_obj}. Partial-bids and offers for storage charging and discharging are accepted in constraints \eqref{eq:lo_chs}--\eqref{eq:lo_dis}. Other constraints on market clearing  \eqref{eq:lo_same} are again identical to those described in the previous section.    

This bilevel formulation is nonlinear and non-convex. However, under the assumption of convexity of the lower level, it can be parametrized using the Karush–Kuhn–Tucker optimality conditions. The complementary slackness conditions are further linearized using the Fortuny-Amat and  McCarl transformation \cite{arroyo_2010}. The nonlinear terms in the objective function are also linearized using the strong duality condition. The resulting single-level equivalent is a mixed-integer linear program that can be solved with commercial solvers. Interested readers are referred to \cite{yishen_2016, yishen_ferc, baringo_2016_windoffer, jalal_2016_pscc} for further details.

\section{Case Study}
\begin{figure} [!b]
\vspace{-10pt}
	\includegraphics[width=0.30\textwidth, height=0.35\textheight]{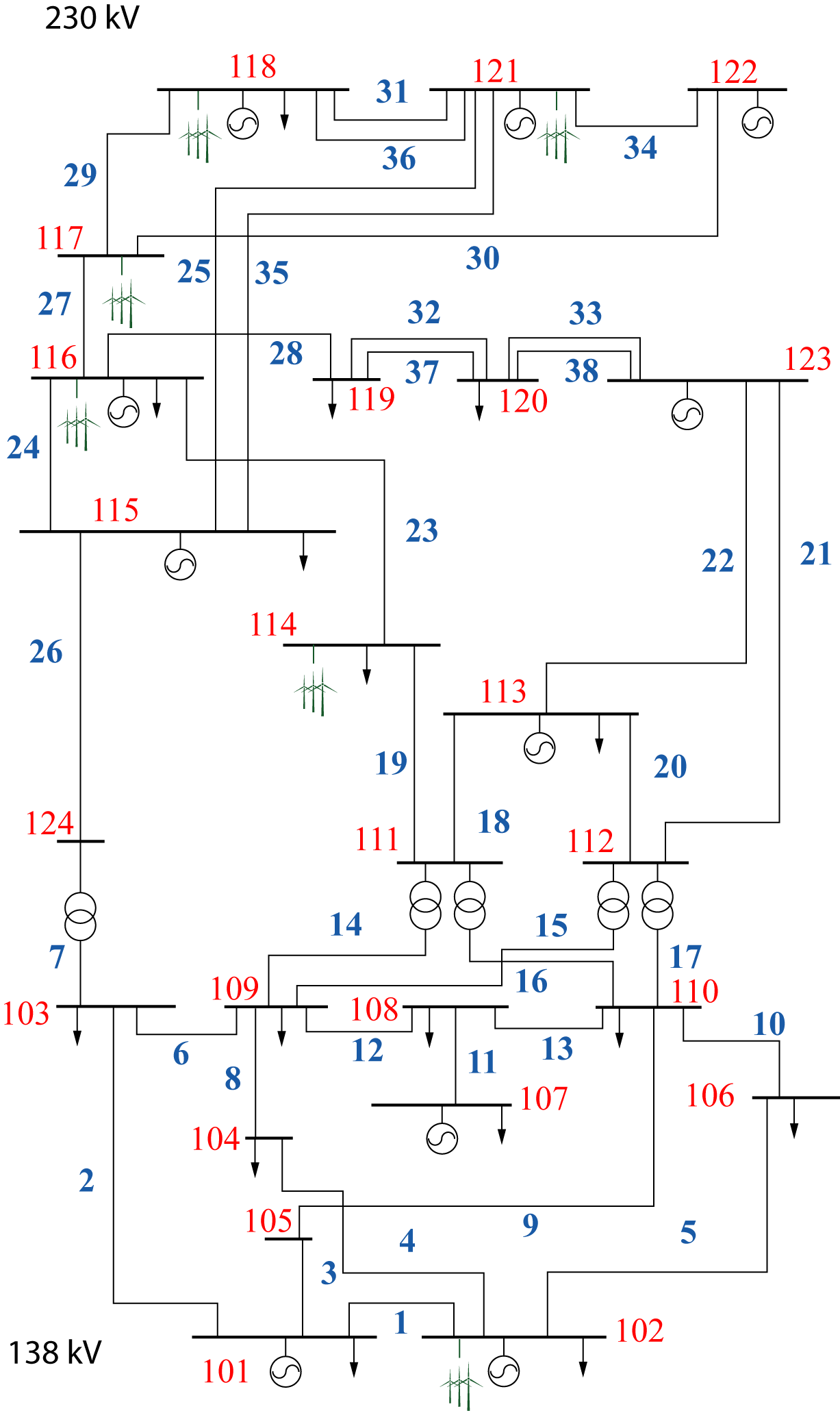}
	\centering
	\caption{IEEE One-Area Reliability Test System.} \label{plot:rts}
\end{figure}

The proposed market models have been tested on the modified version of the IEEE Reliability Test System (RTS) \cite{hrvoje_2015}, Fig.~\ref{plot:rts}. This test system consists of 24 buses, 32 generators, 38 transmission lines, and 5 wind farms. The day-ahead wind power forecasts were generated based on the NREL Eastern Wind Dataset and allow us to assess the proposed models over a one-year simulation horizon. We assume one storage with the optimized parameters from \cite{hrvoje_2015}: $\overline{dis}_h = \overline{chs}_h=93$ MW, $\overline{SoC}_h=629$ MWh, $\eta_h^{dis} = \eta_h^{chs}=0.9$, $SoC_{h}^{init}=315$ MWh, $C_{h}^{dis/chs}=0$ \$/MWh, $C_{h}^{bid}=30$ \$/MWh.

All simulations were carried out using GAMS 23.7 and CPLEX 12.5 on a Intel Xenon 2.55 GHz processor with 32 GB RAM. The computation time for each single day is less than or equal to 5 s with the MILP gap set to 0.005 \%.

To examine the effect of local congestion on the transmission network, we first scale up the original line capacities by 50\% to ensure that there is no congestion in the network, and that the LMP is the same at all buses. We then gradually reduce the capacities of the lines connected to the bus where the storage is located. Other line capacities remain unchanged in order to avoid congestion across the system. We adopted this approach because we wanted to focus on what happens when storage is not able to deliver its flexibility. By comparing operating profits when locating storage at different buses, we show how siting decisions affect storage profitability.


\begin{figure}[!t]
	\includegraphics[width=0.4\textwidth]{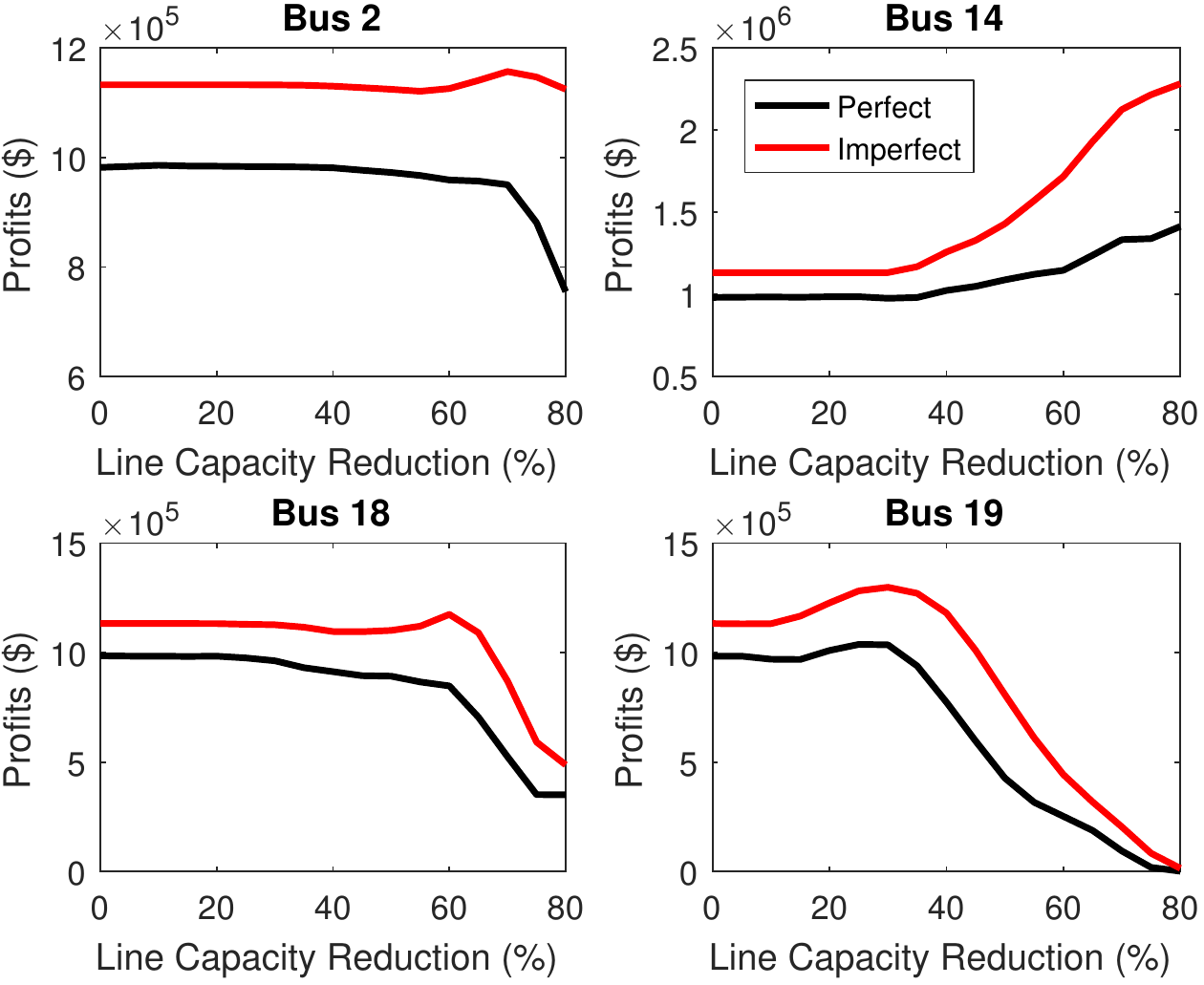}
	\centering
	\caption{Typical patterns of annual storage profits.} 
	\label{plot:line}
		\vspace{-20pt}
\end{figure}

Figure~\ref{plot:line} shows how the annual profit collected by the storage changes as the local line capacities are reduced when this storage is located at four different buses. Decreasing these line capacities tends to increase the LMP, which in turn tends to improve the profitability of storage for all four locations. However, these figures show that different patterns are possible. For example, if the storage is located at buses 2, 18 or 19, at some point the increase in LMP resulting from a reduction in transmission capacity is offset by a reduction in the amount of energy that the storage can physically deliver to the rest of the system. On the other hand, if the storage is located at bus 14, profitability increases monotonically as the local transmission capacity decreases. Factors such as the location of conventional and wind generation, the size and location of the loads, the topology of the network as well as the transmission capacity at nearby buses determine how local congestion affects storage profitability. While patterns for perfect and imperfect competition are similar, strategic bidding significantly increases profitability particularly when local congestion is significant. 

\begin{figure}[!b]
	\includegraphics[width=0.4\textwidth]{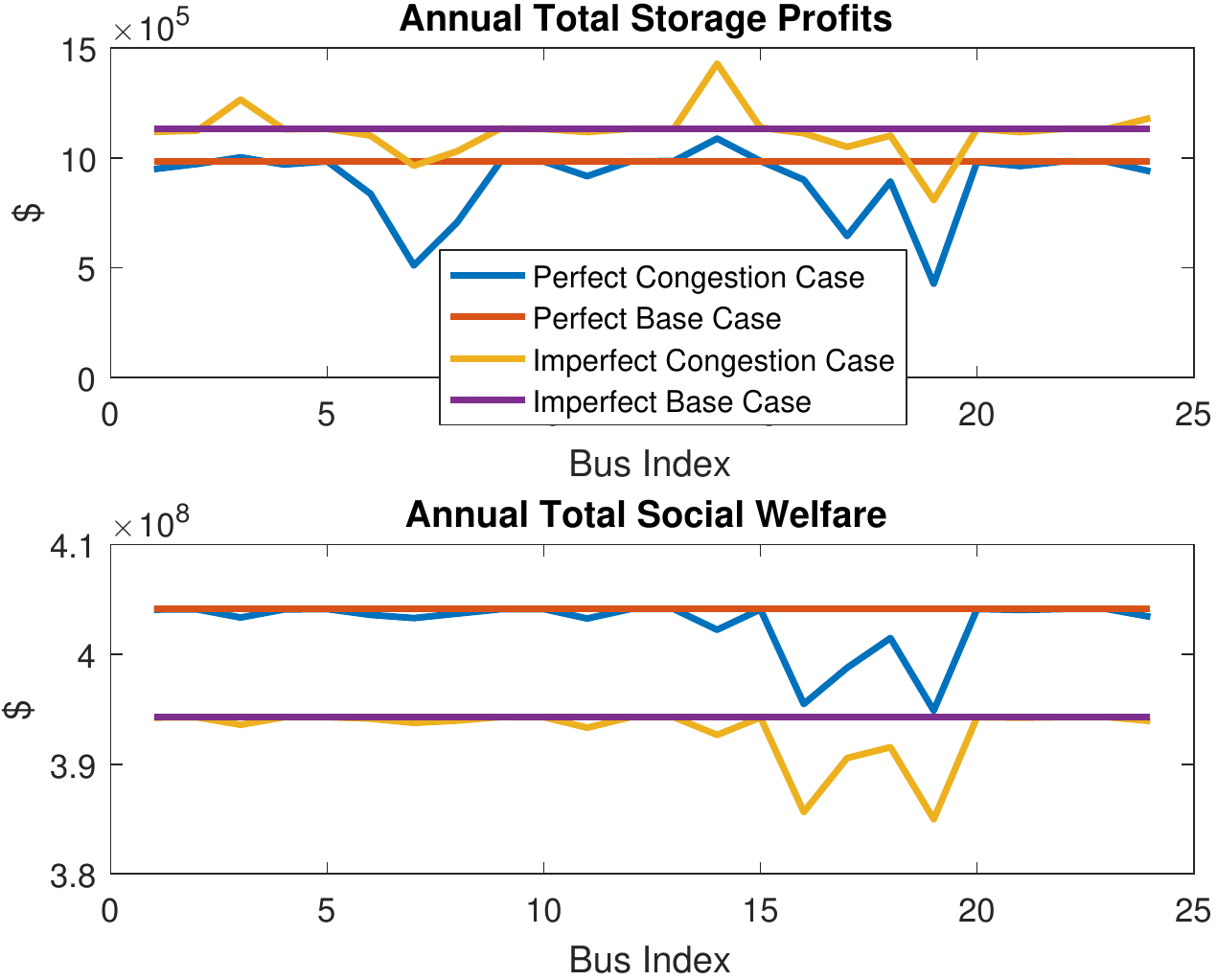}
	\centering
	\caption{Annual profits at various buses under a 50\% reduction in line transmission capacity.} 
	\label{plot:bus}
\end{figure}

\begin{figure} [htb]
	\includegraphics[width=0.4\textwidth]{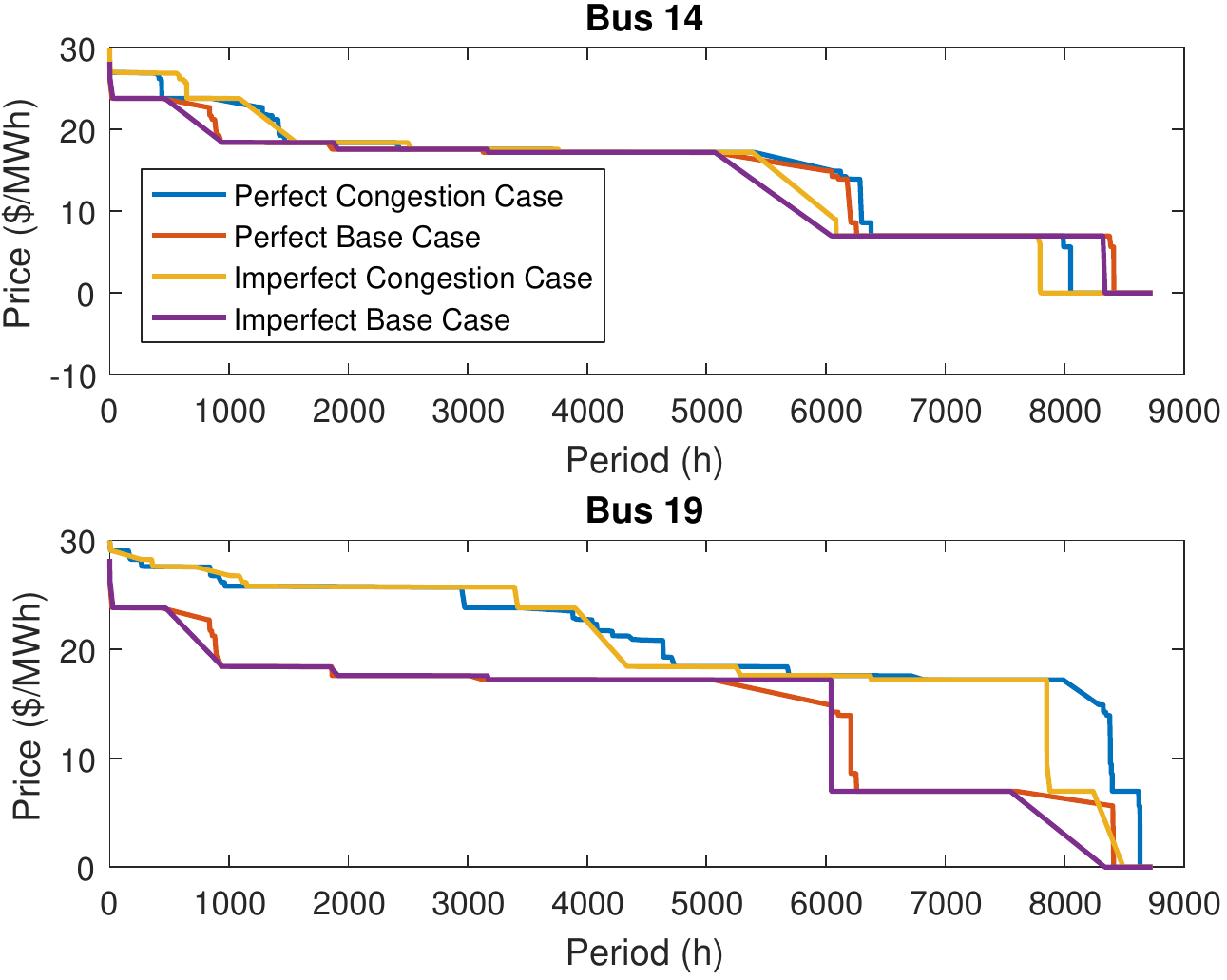}
	\centering
	\caption{Locational marginal price duration curve with a 50\% line capacity reduction.} 
	\label{plot:lmpcdf}
	\vspace{-25pt}
\end{figure}

Figure~\ref{plot:bus} shows the annual profit and social welfare for all possible storage locations under a 50\% line capacity reduction. Imperfect competition enhances storage profitability for all locations but causes a reduction in the social welfare. Compared with the base case, storage can collect more or less profits depending on the location. This pattern is relatively consistent regardless of whether competition is perfect or imperfect. In terms of social welfare, the shapes of both cases are reasonably similar. For buses 15--19, social welfare is more sensitive to the congestion level. Strategic behavior actually does not result in a significant loss of social welfare. These observations motivate the search for the optimal storage locations because such locations would ensure that storage collects enough revenue while being less affected by congestion in the system.

Figure~\ref{plot:lmpcdf} shows the LMP duration curve at buses 14 and 19 based on a year-long simulation for the base case and for a 50\% reduction on line capacity. The prices in perfect and imperfect competition are relatively close, but slightly more favorable for merchant storage in case of imperfect competition. On the other hand, comparing the prices for the base case shows that local congestion increases the average nodal price. The top left area of both plots suggest that high prices happen more frequently than in the base case, while the bottom right areas show that low prices happen less frequently. Significant penetration of renewables causes  periods of zero or negative prices occur and are made more frequent by local congestion as shown in the plot for bus 14. 
\begin{figure}[!b]
	\vspace{-10pt}
	\includegraphics[width=0.4\textwidth]{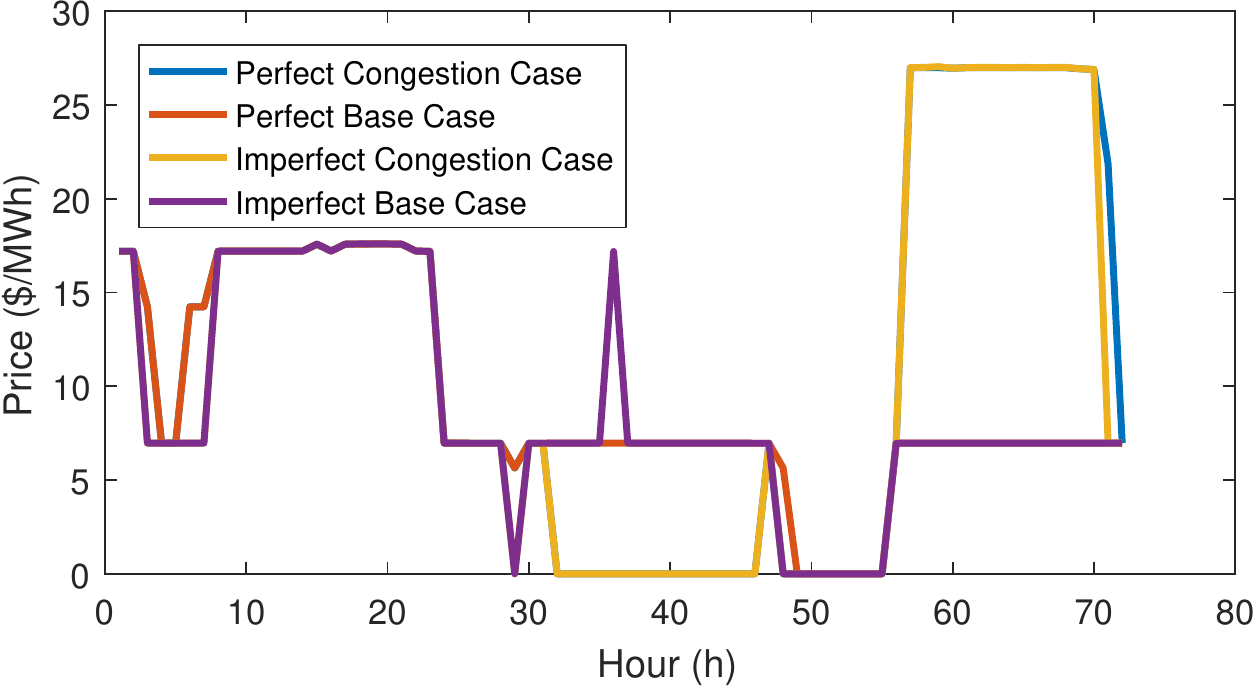}
	\centering
	\caption{LMP at bus 14 on days 100--102 with a 50\% line capacity reduction.} 
	\label{plot:lmpday}
\end{figure}

Figure~\ref{plot:lmpday} shows a three-day snapshot of LMPs at bus 14. During the first 28 hours of this period, local congestion does not affect the LMPs and perfect competition actually produces a slightly higher average LMP. From about hour 28 onwards, strategic bidding can be used to create arbitrage opportunities by driving the LMPs down and up. Local congestion cases enhance these opportunities by extending the periods of lower and higher LMPs. Such price differences are valuable for storage because they make it possible to charge at low prices and discharge at high prices. Nontrivial profit differences between perfectly and imperfectly competitive markets stem from such infrequent periods with substantial price differences. 

\begin{figure} [htb]
\vspace{-5pt}
	\includegraphics[width=0.4\textwidth]{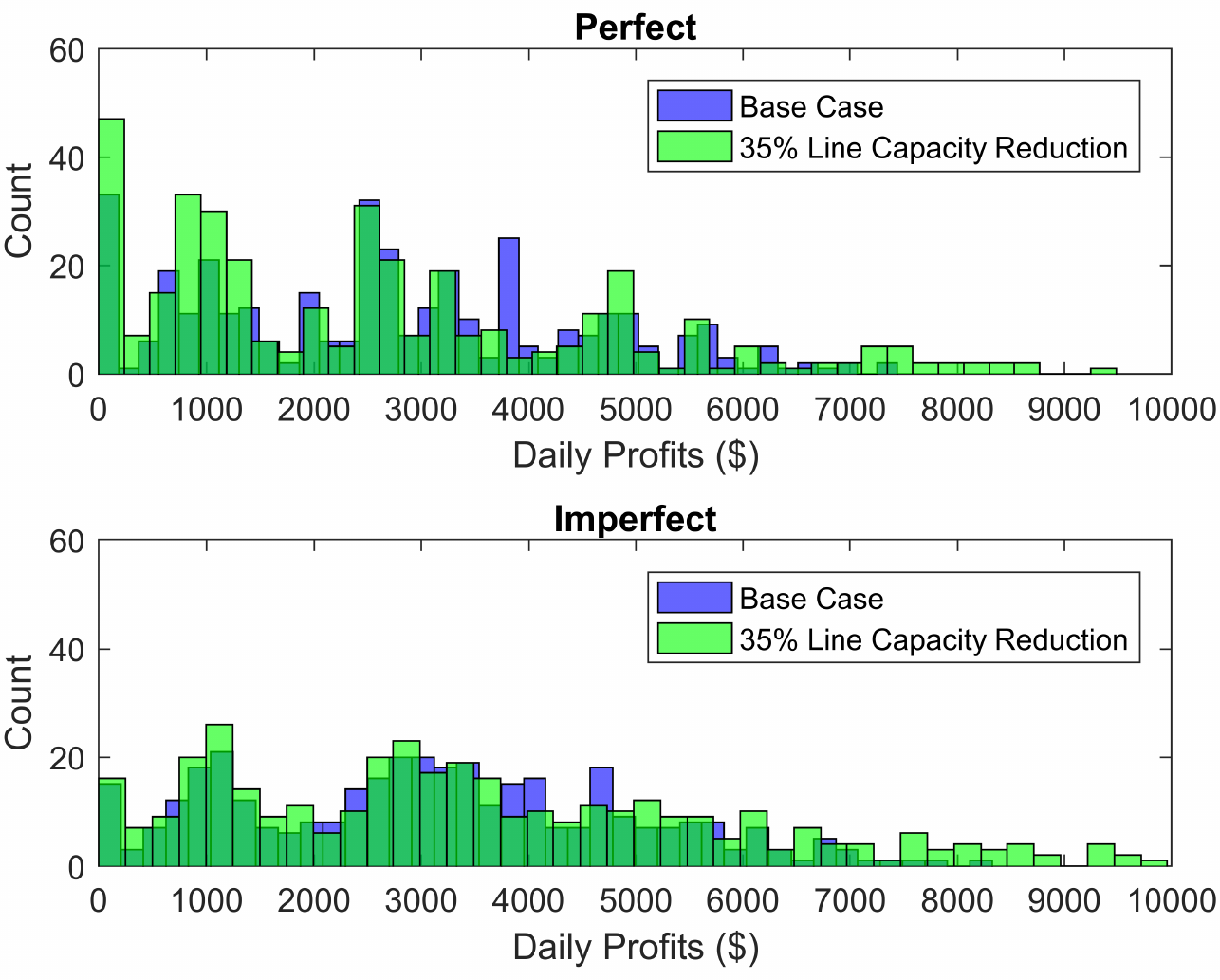}
	\centering
	\caption{Histogram of daily storage profits at bus 19.} 
	\label{plot:profitpdf}
\end{figure}

Figure~\ref{plot:profitpdf} shows the histogram of daily profits for a storage located at bus 19. In the perfect competition case, the histogram does not exhibit a clear pattern. Local congestion increases the number of higher profit days and, to a lesser extent, the number of lower profit days, resulting in a higher total annual profit. The histogram for imperfect competition case is smoother and more evenly distributed. There are much fewer  low profit days and many more high profit days. Here again, congestion increases the number of high profit days. 
 
\vspace{-20pt}
\section{Conclusion}

Local transmission congestion can alternate locational marginal prices substantially. Although these differences can enhance the profitability of energy storage systems,  local congestion can hamper the ability of storage to participate in arbitrage on a transmission system level. This paper quantifies  these effects by simulating the outcomes of market-clearing in a system with a large amount of storage and wind generation under both perfect and imperfect competition. These simulations show that a modest level of local transmission congestion increases nodal market-clearing prices, which results in higher energy storage profits. Depending on the location of the storage, this profitability is sustained even for severe local congestion. Imperfect competition is also shown to help storage increase their profits by influencing the LMPs. The results also show the importance of the storage location on its profitability. 
\vspace{-10pt}
\bibliographystyle{IEEEtran}	
\bibliography{literature}		
	
\end{document}